\documentclass[preprint,3p,times]{elsarticle}

\usepackage{lineno,hyperref}
\usepackage{amsmath}
\usepackage{amssymb}
\usepackage{xcolor}
\usepackage{bbm}
\usepackage{mathtools}
\usepackage{amsthm}
\usepackage{pgfplots}
\usepackage{import}
\usepackage{float}
\usepackage{tikz}
\pgfplotsset{compat = 1.3}

\biboptions{numbers,sort&compress}
\usetikzlibrary{patterns,calc,decorations.pathmorphing,decorations.markings}

\journal{}

\newcommand{\diff}{\mathrm{d}}
\newtheorem{theorem}{Theorem}

\newtheorem{remark}[theorem]{Remark}
\newtheorem{definition}[theorem]{Definition}

\begin{document}
\begin{frontmatter}

\title{Operator splitting for coupled linear port-Hamiltonian systems}

\author[BUW]{Jan Lorenz}
\ead{jan.lorenz@uni-wuppertal.de}

\author[BUW]{Tom Zwerschke}
\ead{tom.zwerschke@uni-wuppertal.de}

\author[BUW]{Michael Günther}
\ead{guenther@uni-wuppertal.de}

\author[BUW]{Kevin Schäfers\corref{cor1}}
\ead{kschaefers@uni-wuppertal.de}
\cortext[cor1]{Corresponding author.}

\affiliation[BUW]{organization={Bergische Universität Wuppertal, IMACM},
             addressline={Gaußstraße 20},
             city={Wuppertal},
             postcode={D-42119},
             country={Germany}}

\begin{abstract}
Operator splitting methods tailored to coupled linear port-Hamiltonian systems are developed. We present algorithms that are able to exploit scalar coupling, as well as multirate potential of these coupled systems. The obtained algorithms preserve the dissipative structure of the overall system and are convergent of second order. Numerical results for coupled mass-spring-damper chains illustrate the computational efficiency of the splitting methods compared to a straight-forward application of the implicit midpoint rule to the overall system.
\end{abstract}


\end{frontmatter}

\section{Introduction}\label{sec:introduction} 
\noindent Many (complex) physical systems can be modeled as linear port-Hamiltonian systems (PHS) of ordinary differential equations (ODEs) \cite{van2006port}
\begin{subequations}\label{eq:linear_PHS}
\begin{align}
    \dot{x}(t) &= \left[J -R\right]\nabla H(x) + B u(t), \quad x(t_0) = x_0, \\
    y(t) &= B^\top \nabla H(x),
\end{align}
\end{subequations}
with state variable $x \colon [t_0,t_\mathrm{end}] \to \mathbb{R}^n$, quadratic Hamiltonian $ H(x) = \tfrac{1}{2} x^\top Q x$, $J = -J^\top \in \mathbb{R}^{n \times n}$ skew-symmetric, $R \in \mathbb{R}^{n \times n}$ positive semi-definite, input $u \colon [t_0,t_\mathrm{end}] \to \mathbb{R}^d$ and output $y \colon [t_0,t_\mathrm{end}] \to \mathbb{R}^d$, as well as the port matrix $B \in \mathbb{R}^{n \times d}$. 
Supposing sufficient regularity of the input, there exists an unique solution that satisfies the energy balance condition
\begin{align} \label{eq:energy_balance}
    \tfrac{\diff}{\diff t} H(x(t)) = -\nabla H(x(t))^\top R(x(t)) \nabla H(x(t)) + y(t)^\top u(t) \leq y(t)^\top u(t).
\end{align}
The dissipation inequality \eqref{eq:energy_balance} can also be stated in integral form, i.e.,
$ H(x(t_0 + h)) - H(x(t_0)) \leq \int\nolimits_{t_0}^{t_0 + h} y(\tau)^\top u(\tau)\;\diff \tau$.
We demand numerical solutions of \eqref{eq:linear_PHS} that preserve the energy balance \eqref{eq:energy_balance} at a discrete level.
This can be achieved, for example, by using discrete gradient methods \cite{gonzalez1996time}.
Furthermore, we are aiming at deriving efficient numerical integration schemes which exploit the fine structure of the PHS.
In this paper, we consider systems \eqref{eq:linear_PHS} with finer structure 
\begin{subequations}\label{eq:overall-system}
\begin{align}
    \dot{x}(t) &= \left[\!\begin{pmatrix}
        J_{1} & -\tilde{J}^\top \\ \tilde{J} & J_2 
    \end{pmatrix} - \begin{pmatrix}
        R_1 & 0 \\ 0 & R_2
    \end{pmatrix}\!\right] \begin{pmatrix}
        Q_1 & 0 \\ 0 & Q_2
    \end{pmatrix} x + Bu(t), \quad x(t_0) = x_0, \\
    y(t) &= B^\top \nabla H(x),
\end{align}
\end{subequations}
where $J_1 \in \mathbb{R}^{n_1 \times n_1},\; J_2 \in \mathbb{R}^{n_2 \times n_2}$ 
are skew-symmetric, $R_1 \in \mathbb{R}^{n_1 \times n_1},\; R_2 \in \mathbb{R}^{n_2 \times n_2}$ are positive semi-definite matrices, and $Q_1 \in \mathbb{R}^{n_1 \times n_1},\; Q_2 \in \mathbb{R}^{n_2 \times n_2}$, $n = n_1 + n_2$. 
This structure typically arises by coupling linear PHS via port matrices, e.g.\ in the modeling of electric circuits~\cite[Remark 7]{bartel2023operator}, and then shifting the coupling in the off-diagonal blocks of the skew-symmetric matrix.

Operator splitting techniques \cite{mclachlan2002splitting} enable the exploitation of the particular structure by additively splitting the right-hand side of \eqref{eq:overall-system}.
Recently, operator splitting techniques have been applied successfully to linear port-Hamiltonian ODE systems by splitting the PHS into its energy-conservative and energy-dissipative parts \cite{frommer2023operator}.
In this work, we will consider different splitting approaches that split the overall system \eqref{eq:overall-system} into its subsystems.
Nonetheless, we want to emphasize that the resulting subsystems in our approach can be split further into their conservative and dissipative parts by applying the approaches introduced in \cite{frommer2023operator} in a hierarchical manner.

The paper is organized as follows. 
In Section \ref{sec:operator_splitting} we give a brief overview on operator splitting techniques and dissipativity-preserving numerical integration schemes for linear PHS.
Section \ref{sec:scalar_coupling} is devoted to the special case where the matrix $\tilde{J}$ introduces a scalar coupling. 
Often, the port-Hamiltonian subsystems have transient behavior that is characterized by different time scales. In Section \ref{sec:multirate}, we will develop multiple time stepping schemes that are integrating the subsystems using different step sizes.
The proposed splitting methods are applied to coupled mass-spring-damper chains in Section \ref{sec:numerical_results}.
The paper closes with a summary and an outlook for future research.

\section{Operator Splitting for Linear PHS}\label{sec:operator_splitting}
\noindent Operator splitting \cite{mclachlan2002splitting} is a framework to efficiently solve additively partitioned initial-value problems (IVPs) of ODEs 
\begin{align}\label{eq:split-system}
    \dot{x}(t) = f(x,t) = f^{\{1\}}(x,t) + f^{\{2\}}(x,t), \quad x(t_0) = x_0,
\end{align}
where the right-hand side is split into two different parts w.r.t., for example, stiffness, nonlinearity, dynamical behavior and computational cost. 
Rewriting \eqref{eq:split-system} in the homogeneous form 
\begin{equation}\label{eq:dynamical_system}
    \begin{pmatrix}
        \dot{x} \\ \dot{s} 
    \end{pmatrix} = \begin{pmatrix}
        f(x,s) \\ 1
    \end{pmatrix} = \underbrace{\begin{pmatrix}
        f^{\{1\}}(x,s) \\ 1
    \end{pmatrix}}_{\eqqcolon \tilde{f}^{\{1\}}(x,s)} + \underbrace{\begin{pmatrix}
        f^{\{2\}}(x,s) \\ 0
    \end{pmatrix}}_{\eqqcolon f^{\{2\}}(x,s)}, \quad \begin{pmatrix}
        x(t_0) \\ s(t_0)
    \end{pmatrix} = \begin{pmatrix}
        x_0 \\ 0
    \end{pmatrix},
\end{equation}
operator splitting techniques approximate the exact flow $\varphi_t(x_0)$ of the dynamical system \eqref{eq:dynamical_system} by composing the flows $\varphi_t^{\{1\}}(x_0)$ and $\varphi_t^{\{2\}}(x_0)$ of the subsystems $ \dot{x} = \tilde{f}^{\{1\}}(x,s) $ and  $ \dot{x} = \tilde{f}^{\{2\}}(x,s) $, respectively \cite{frommer2023operator}. 
Splitting methods of order $p>2$ require negative weights \cite{suzuki1991general}, which violates the dissipativity inequality and thus introduces an order barrier. The paper aims at developing efficient splitting methods of order $p=2$ that are tailored to the specific structure of the linear PHS \eqref{eq:overall-system}. 
For this purpose, we consider the Strang splitting \cite{strang1968construction}
\begin{align}\label{eq:Strang-splitting}
    \Psi_h = \varphi_{h/2}^{\{1\}} \circ \varphi_h^{\{2\}} \circ \varphi_{h/2}^{\{1\}}.
\end{align}
Assuming that both subsystems are linear PHS, splitting methods satisfy the dissipativity inequality at a discrete level since the composition of dissipativity preserving maps is again dissipativity preserving. 

Instead of using the exact flows $\varphi_t^{\{1\}}(x_0)$ and $\varphi_t^{\{2\}}(x_0)$, we can replace them by numerical approximations $\Phi_h^{\{1\}}(x_0)$ and $\Phi_h^{\{2\}}(x_0)$, respectively.
If these are obtained by discrete gradient methods of order $p=2$, it does not affect the convergence order, as well as the energy balance of the splitting method.
For linear PHS \eqref{eq:linear_PHS}, the implicit midpoint rule defines a discrete gradient method of order $p=2$,
\begin{subequations}
\begin{align}\label{eq:implicit-midpoint}
\begin{split}
    x_1 &= x_0 + h [J - R] Q \tfrac{x_0 + x_1}{2} + h B u_1, \\
    y_1 &= B^\top Q\tfrac{x_0 + x_1}{2}.
\end{split}    
\end{align}
with $u_1 = u(\tfrac{t_0 + t_1}{2})$.
In particular, the midpoint $Q \tfrac{x_0+x_1}{2}$ is a discrete approximation of the gradient $\nabla H(x) = Qx$ for the quadratic Hamiltonian $H(x) = \tfrac{1}{2} x^\top Q x$. As a discrete gradient, it satisfies $(x_1 - x_0)^\top Q \tfrac{x_0 + x_1}{2} = H(x_1) - H(x_0)$, which implies the energy balance~\cite{Maccelli2023} 
\begin{align*}
    \frac{H(x_1)-H(x_0)}{h} = \tfrac{(x_1-x_0)^\top}{h} Q \tfrac{x_1+x_0}{2} = \left([J - R] Q \tfrac{x_0 + x_1}{2} + B u_1\right)^\top Q \tfrac{x_1+x_0}{2} = -\left(\tfrac{x_0+x_1}{2}\right)^\top Q^\top R^\top Q \tfrac{x_0+x_1}{2} + u_1^\top y_1 \le u_1^\top y_1.
\end{align*}
By defining $A \coloneqq (J-R)Q$, the update of the state variable can be written as
\begin{align}\label{eq:midpoint4linear}
    \left(I-\tfrac{h}{2}A\right) x_1 &= \left(I+ \tfrac{h}{2}A\right) x_0 + h B u_1 & &\Leftrightarrow & x_1 &= \left(I-\tfrac{h}{2}A\right)^{-1} \left(I+ \tfrac{h}{2}A\right) x_0 + h \left(I - \tfrac{h}{2}A\right)^{-1} B u_1,
\end{align}
\end{subequations}
emphasizing its close connection to the Cayley transform. 
The next two sections will discuss different splittings of the partitioned 
linear PHS \eqref{eq:overall-system} exploiting the fine structure of the system, resulting in an efficient computational process.

\section{Coupled Linear PHS with Scalar Coupling}\label{sec:scalar_coupling}
\noindent In this section, we consider the linear PHS \eqref{eq:overall-system} as an additively partitioned ODE system \eqref{eq:split-system} with splitting
\begin{align} \label{eq:diag_splitting}
    f^{\{1\}}(x,t) &= \underbrace{\begin{pmatrix}
        0 & - \tilde{J}^\top \\ \tilde{J} & 0
    \end{pmatrix}}_{A_1} \nabla H(x) + Bu(t), & f^{\{2\}}(x,t) &= \underbrace{\begin{pmatrix}
        J_1 - R_1 & 0 \\ 0 & J_2 - R_2
    \end{pmatrix}}_{A_2} \nabla H(x).
\end{align}
This splitting preserves the port-Hamiltonian structure in each subsystem, allowing the preservation of the energy balance \eqref{eq:energy_balance} of the system using appropriate numerical integration schemes. Moreover, we can consider the subsystem $\dot{x} = f^{\{2\}}(x,t)$ 
as the internal evolution equation of the systems and $\dot{x} = f^{\{1\}}(x,t)$ 
as the coupling part.\newline
\noindent\textbf{Computational complexity.}
Considering the implicit midpoint rule and the computation of the update \eqref{eq:midpoint4linear} via an LU-decomposition, the splitting approach is generally more complex than the straight-forward application of the implicit midpoint rule to the overall system. 
The situation changes if we make further assumptions on the coupling.

\begin{definition}
The system \eqref{eq:diag_splitting} is scalar coupled, if the coupling matrix $A_1 Q$ has just one non-zero above, as well as below the diagonal.
\end{definition}
Scalar coupling occurs in many use cases of linear PHS. 
We can exploit this structure to reduce the computational cost to compute \eqref{eq:midpoint4linear} for the second subsystem. 
As we have seen in \eqref{eq:midpoint4linear}, the implicit midpoint rule demands the solution of a linear system of equations. Scalar coupling enables us to compute the transformation directly in a numerical stable way with no need of a matrix decomposition.
Let ${\tilde j}_1$ be the non-zero value below and $\tilde{j}_2$ the non-zero value above the diagonal of $A_1 Q$.
Then, it holds
\begin{align*}
 &\begin{pmatrix} 1 & & & \hdots & & & 0 \\
		 & \ddots & & & & &  \\
                 & & 1 & \hdots & -{\tilde j}_2 & & \\
		 \vdots & & \vdots & \ddots & \vdots & & \vdots \\
                 & & -{\tilde j}_1 & \hdots & 1 & & \\
		 & & & & & \ddots & \\
	0 & & & \hdots & & & 1 \end{pmatrix} ^ {-1} 
 \begin{pmatrix} 1 & & & \hdots & & & 0 \\
		 & \ddots & & & & &  \\
                 & & 1 & \hdots & {\tilde j}_2 & & \\
		 \vdots & & \vdots & \ddots & \vdots & & \vdots \\
                 & & {\tilde j}_1 & \hdots & 1 & & \\
		 & & & & & \ddots & \\
	0 & & & \hdots & & & 1 \end{pmatrix} 
 = &\begin{pmatrix} 1 & & & \hdots & & & 0 \\
		 & \ddots & & & & &  \\
		 & & \frac{1 + \tilde{j}_1\tilde{j}_2}{1 - \tilde{j}_1 \tilde{j}_2} & \hdots & \frac{2\tilde{j}_2}{1 - \tilde{j}_1 \tilde{j}_2} & & \\
		 \vdots & & \vdots & \ddots & \vdots & & \vdots \\
		 & & \frac{2 \tilde{j}_1}{1 - \tilde{j}_1 \tilde{j}_2} & \hdots & \frac{1 + \tilde{j}_1 \tilde{j}_2}{1 - \tilde{j}_1 \tilde{j}_2} & & \\
		 & & & & & \ddots & \\
	0 & & & \hdots & & & 1 \end{pmatrix}.
\end{align*}
Applying the implicit midpoint rule \eqref{eq:implicit-midpoint} to the overall system \eqref{eq:linear_PHS} demands the computation of $(I - \tfrac{h}{2} A)$ and $(I + \tfrac{h}{2} A)$, and the LU decomposition of the first matrix, resulting in setup costs of $\mathcal{O}(n^3)$. Moreover, one time step of the implicit midpoint rule requires $\mathcal{O}(n^2 + nd)$ operations. 
In contrast, the splitting method \eqref{eq:Strang-splitting} with splitting \eqref{eq:diag_splitting} and assuming scalar coupling has setup costs of size $\mathcal{O}(n_1^3 + n_2^3)$ and demands $\mathcal{O}(nd + n_1^2 + n_2^2)$ operations per time step.
Since it usually holds $d \ll n$, the splitting approach becomes more efficient. Particularly, a straight-forward computation shows that the splitting method requires less operations, provided that $d \le n$ and $n_1,n_2 \ge 2$.

\section{Coupled Linear PHS with Multirate Potential}\label{sec:multirate}
\noindent In this section, we consider linear PHS \eqref{eq:overall-system}, and we introduce the splitting approach
\begin{align}\label{eq:multirate-splitting}
    \dot{x} &= f^{\{1\}}(x,t) + f^{\{2\}}(x,t) &
    f^{\{1\}}(x,t) &=  \underbrace{
    \begin{pmatrix}
        0 & -\tilde{J}^\top \\
        \tilde J & J_2 - R_2
    \end{pmatrix}}_{M_1}   
    \nabla H(x) + Bu(t), \quad
    f^{\{2\}}(x,t) = \underbrace{
    \begin{pmatrix}
        J_1-R_1 & 0 \\
        0 & 0
    \end{pmatrix}}_{M_2} 
    \nabla H(x),
\end{align}
where $f^{\{1\}}$ is characterized by slow dynamics and $f^{\{2\}}$ by fast dynamics.
Moreover, we assume $n_1 \ll n_2$. 
To exploit this sort of multirate behavior, we aim at using a large macro-step size $h$ for the slow subsystem 
\begin{subequations}
    \begin{equation}\label{eq:slow-subsystem}
        \dot{x} = M_1 \nabla H(x) + Bu(t),
    \end{equation}
and a small micro-step size $h/m,\ m\in\mathbb{N},\ m \gg 1,$ for the fast subsystem
\begin{equation}\label{eq:fast-subsystem}
    \dot{x} = M_2 \nabla H(x). 
\end{equation}
\end{subequations}
This can be achieved by using the \emph{impulse method}, a multirate generalization to the Strang splitting \cite{HairerLubichWanner}.
It reads
\begin{equation}\label{eq:impulse_method}
    \Psi_h = \varphi_{h/2}^{\{1\}} \circ \left( \varphi_{h/m}^{\{2\}}\right)^m \circ \varphi_{h/2}^{\{1\}}
\end{equation}
and computes the numerical approximation $x_1 \approx x(h)$ by composing (approximated) flows of the subsystems \eqref{eq:slow-subsystem} and \eqref{eq:fast-subsystem}. 
Since the integrator \eqref{eq:impulse_method} is a composition of flows of the subsystems \eqref{eq:slow-subsystem} and \eqref{eq:fast-subsystem} that are port-Hamiltonian and all weights are positive, the numerical approximation satisfies the dissipativity inequality.\newline

\noindent\textbf{Computational complexity.} 
Due to the simple structure of $M_2$, the fast subsystem \eqref{eq:fast-subsystem} simplifies to 
\begin{equation}
    \dot{x}_1 = \tilde{M}_2 \nabla_{1} H(x),
\end{equation} 
with $\tilde{M}_2 \coloneqq J_1 - R_1$, $x_1$ denoting the first $n_1$ components and $\nabla_{1} H(x)$ denoting the first $n_1$ components of $\nabla H(x)$. 

Neglecting the setup costs, computing one time step of the impulse method \eqref{eq:impulse_method} with macro-step size $h$ demands $\mathcal{O}(n^2 + m n_1^2)$ operations. 
In contrast, computing $m$ steps of the implicit midpoint rule \eqref{eq:implicit-midpoint} with step size $h/m$ to resolve the fast dynamics accurately enough demands $\mathcal{O}(mn^2)$ operations. 
If the system \eqref{eq:multirate-splitting} exhibits multirate potential so that $n_1 \ll n_2$ and $m \gg 1$, the multiple time stepping approach requires less operations.

\remark{An alternative approach would shift the coupling part into the fast subsystem. However, this would imply that the fast subsystem is of full dimension $n = n_1 + n_2$ so that the $m \gg 1$ evaluations of the fast flow $\varphi^{\{2\}}$ would become more expensive, resulting in a less efficient computational process.} \normalfont

\section{Numerical Results}\label{sec:numerical_results}
\noindent A frequently used benchmark problem for linear PHS is given by a mass-spring-damper chain\footnote{\href{https://algopaul.github.io/PortHamiltonianBenchmarkSystems.jl/}{https://algopaul.github.io/PortHamiltonianBenchmarkSystems.jl/}} (MSD-chain)~\cite{GUGERCIN20121963}.
We consider two coupled MSD-chains, consisting of $n_1$ and $n_2$ mass-spring-damper elements, respectively. The two chains are coupled via a spring with stiffness $K_{co} > 0$, see Fig.\ \ref{fig:coupled_MSD-chain}. 
The first chain consists of homogeneous masses $m_1 > 0$, mass-less springs $K_1 > 0$, and damping $r_1 > 0$. For the second chain, we have homogeneous parameters $m_2 > 0$, $K_2 > 0$, and $r_2 > 0$.
Introducing positions $q_{i,1},\ldots,q_{i,n_i}$ and momenta $p_{i,1},\ldots,p_{i,n_i}$, $i=1,2$, we can formulate the equations of motion as a PHS of the form \eqref{eq:overall-system} by setting $q = q_{1,1}$ as the coupling variables, resulting in \cite{bartel2023operator}
\begin{align*}
    J_1 &= \begin{pmatrix}
    0 & -1 &  & & & -1 \\
    1 & 0 &  & & & 0 \\
    & & \ddots & & & \vdots \\
    & &  & 0 & -1 &  \\
    & &  & 1 & 0 & 0 \\
    1 & 0 & \cdots &  & 0 & 0
\end{pmatrix}\in \mathbb{R}^{(2n_1+1) \times (2n_1+1)}, & J_2 &= \begin{pmatrix}
    0 & -1 &  & & &  \\
    1 & 0 &  & & &  \\
    & & & \ddots & &   \\
    & & &  & 0 & -1   \\
    & & &  & 1 & 0 \\
\end{pmatrix}\in \mathbb{R}^{2n_2 \times 2n_2}, \\
R_1 &= \mathrm{diag}(r_{1},0,r_{1},0,\ldots,r_{1},0,0)\in \mathbb{R}^{(2n_1+1) \times (2n_1+1)}, & R_2 &= \mathrm{diag}(r_{2},0,r_{2},0,\ldots,r_{2},0)\in \mathbb{R}^{2n_2 \times 2n_2}, \\
\tilde{J} &= \begin{pmatrix}
    0 & \cdots & 0 & 1 \\
    0 & \cdots & 0 & 0 \\
    \vdots & & \vdots & \vdots \\
    0 & \cdots & 0 & 0
\end{pmatrix}\in \mathbb{R}^{2n_2 \times (2n_1+1)}, & Q &= \begin{pmatrix}
    Q_1 & 0 & 0 \\
    0 & K_{co} & 0 \\
    0 & 0 & Q_2
\end{pmatrix} \in \mathbb{R}^{(2n+1)\times (2n+1)}, 
\end{align*} 
where the two matrices $Q_1,Q_2$ are given by 

$$Q_i = \begin{pmatrix}
    \tfrac{1}{m_{i}} & 0 & 0 & & & & & & \\
    0 & K_{i} & 0 & -K_{i} & & & & & \\
    0 & 0 & \tfrac{1}{m_{i}} & 0 & 0 & & & & \\
    & -K_{i} & 0 & 2K_{i}  & 0 & -K_{i} & & & \\
    & & \ddots & & \ddots & & \ddots & & \\
    & & & \ddots & & \ddots & & 0 &  \\
    & & & & \ddots & & 2K_i & 0 & -K_{i} \\
    & & & & & 0 & 0 & \tfrac{1}{m_{i}} & 0 \\
    & & & & &  & -K_{i} & 0 & 2K_{i}\\
\end{pmatrix} \in \mathbb{R}^{2n_i \times 2n_i}, \quad i=1,2.$$
The state variables are $x = (p_{1,1},q_{1,1},\ldots,p_{1,n_1},q_{1,n_1},q_{1,1}-q,p_{2,1},q_{2,1},\ldots,p_{2,n_2},q_{2,n_2})^\top \in \mathbb{R}^{2n+1}$.

\begin{figure}[htb]
    \centering
    \scalebox{.75}{ \begin{tikzpicture}

 \tikzstyle{spring}=[thick,decorate,decoration={zigzag,pre length=0.4cm,post
 length=0.4cm,segment length=4}]

 \tikzstyle{damper}=[thick,decoration={markings,  
   mark connection node=dmp,
   mark=at position 0.5 with 
   {
     \node (dmp) [thick,inner sep=0pt,transform shape,rotate=-90,minimum
 width=10pt,minimum height=2pt,draw=none] {};
     \draw [thick] ($(dmp.north east)+(2pt,0)$) -- (dmp.south east) -- (dmp.south west) -- ($(dmp.north west)+(2pt,0)$);
     \draw [thick] ($(dmp.north)+(0,-3pt)$) -- ($(dmp.north)+(0,3pt)$);
   }
 }, decorate]

 \tikzstyle{ground}=[fill,pattern=north east lines,draw=none,minimum
		 width=0.55cm,minimum height=0.25cm]
 		\node (wall) [ground, rotate=-90, minimum width=3cm,yshift=-3cm,xshift=0.cm] {};
 \draw (wall.north east) -- (wall.north west);
 
 \node[draw,outer sep=0pt,thick] (MI25) at (-0.6,0)[minimum width=1.25cm, minimum height=1.25cm] {$m_{1}$};
 \node[draw,outer sep=0pt,thick] (MI24) at (2.4,0) [minimum width=1.25cm, minimum height=1.25cm] {$m_{1}$};
 \draw[spring] ($(MI25.east) - (0,0.35)$) -- ($(MI24.west) - (0,0.35)$)  node [midway,below,yshift=-4.5] {$K_{1}$};

 \node (wallI24) [ground, rotate=-90, minimum width=.5cm,yshift=.65cm,xshift=-0.35cm] {};
\draw (wallI24.north east) -- (wallI24.north west);

 \draw[damper] ($(wallI24.north)$) -- ($(MI24.west) + (0,0.35)$)   node  [midway,above,yshift=7.5] {$r_{1}$};
 \draw[thick, dashed] ($(MI25.south)$) -- ($(MI25.south) - (0,0.5)$);
 \draw[thick, dashed] ($(MI24.south)$) -- ($(MI24.south) - (0,0.5)$);
 \draw[thick, -latex] ($(MI25.south) - (0,0.49)$) -- 
                            ($(MI25.south) - (1,0.49)$)
                            node [midway, below] {$q_{1,n_1}$};
 \draw[thick, -latex] ($(MI24.south) - (0,0.49)$) -- 
							($(MI24.south) - (1,0.49)$)
							node [midway, below] {$q_{1,n_1-1}$};
 \draw[spring] ($(MI25.west) - (0,0.35)$) -- ($(wall.north)-(0,.35)$)
 							node [midway,below,yshift=-4.5] {$K_{1}$}; 

 \draw[damper] ($(wall.north)+(0,.35)$) -- ($(MI25.west) + (0,0.35)$)  
							node [midway,above,yshift=7.5] {$r_{1}$}; 

\node at (3.6,-.35) {$\cdots$};


\node[draw,outer sep=0pt,thick] (MI1) at (6.3,0) [minimum width=1.25cm, minimum height=1.25cm] {$m_{1}$};
 \draw[spring] ($(MI1.west) - (0,0.35)$) -- ($(MI1.west)-(1.6,.35)$)
node [midway,below,yshift=-4.5] {$K_{1}$};

\node (wallI1) [ground, rotate=-90, minimum width=.5cm,yshift=4.6cm,xshift=-0.35cm] {};
\draw (wallI1.north east) -- (wallI1.north west);
 \draw[damper] ($(wallI1.north)$) -- ($(MI1.west) + (0,0.35)$)   node  [midway,above,yshift=7.5] {$r_{1}$};

 \draw[thick, dashed] ($(MI1.south)$) -- ($(MI1.south) - (0,0.5)$);
 \draw[thick, -latex] ($(MI1.south) - (0,0.49)$) -- ($(MI1.south) - (1,0.49)$)
node [midway, below] {$q_{1,1}$};


\node[draw,outer sep=0pt,thick] (MII1) at (9.2,0) [minimum width=1.25cm, minimum height=1.25cm] {$m_{2}$};
 \draw[spring] ($(MI1.east) - (0,0.35)$) -- ($(MII1.west)-(0,.35)$)
node [midway,below,yshift=-4.5] {$K_{co}$}; 

 \draw[spring] ($(MII1.east) - (0,0.35)$) -- ($(MII1.east)-(-1.6,.35)$)
node [midway,below,yshift=-4.5] {$K_{2}$}; 

\node (wallII1) [ground, rotate=90, minimum width=0.5cm, 
				yshift=-11.cm,xshift=0.35cm] {};
				\draw (wallII1.north east) -- (wallII1.north west);
 \draw[damper] ($(wallII1.north)$) -- ($(MII1.east) + (0,0.35)$)   node  [midway,above,yshift=7.5] {$r_{2}$};

\node at (11.9,-.35) {$\cdots$};


 \node[draw,outer sep=0pt,thick] (MII25) at (16.1,0)[minimum width=1.25cm, minimum height=1.25cm] {$m_{2}$};
  \node[draw,outer sep=0pt,thick] (MII24) at (13.,0) [minimum width=1.25cm, minimum height=1.25cm] {$m_{2}$};
 \draw[spring] ($(MII24.east) - (0,0.35)$) -- ($(MII25.west) - (0,0.35)$)  node [midway,below,yshift=-4.5] {$K_{2}$};
 
\draw[thick, dashed] ($(MII25.south)$) -- ($(MII25.south) - (0,0.5)$);
\draw[thick, dashed] ($(MII24.south)$) -- ($(MII24.south) - (0,0.5)$);
\draw[thick, -latex] ($(MII25.south) - (0,0.49)$) -- ($(MII25.south) - (-1,0.49)$)
node [midway, below] {$q_{2,n_2}$};
\draw[thick, -latex] ($(MII24.south) - (0,0.49)$) --($(MII24.south) - (-1,0.49)$) node [midway, below] {$q_{2,n_2-1}$};

\node (wallII24) [ground, rotate=90, minimum width=0.5cm, 
yshift=-14.7cm,xshift=0.35cm] {};
\draw (wallII24.north east) -- (wallII24.north west);
\draw[damper] ($(wallII24.north)$) -- ($(MII24.east) + (0,0.35)$)   node  [midway,above,yshift=7.5] {$r_{2}$};

\draw[thick, dashed] ($(MII1.south)$) -- ($(MII1.south) - (0,0.5)$);
\draw[thick, -latex] ($(MII1.south) - (0,0.49)$) -- ($(MII1.south) - (-1,0.49)$) node [midway, below] {$q_{2,1}$};


\node (wall2) [ground, rotate=90, minimum width=3cm, 
				yshift=-18.5cm,xshift=0.cm] {};
\draw (wall2.north east) -- (wall2.north west);

 \draw[spring] ($(MII25.east) - (0,0.35)$) -- ($(wall2.north)-(0,0.35)$)
							node [midway,below,yshift=-4.5] {$K_{2}$}; 
\draw[damper] ($(wall2.north)+(0,.35)$) -- ($(MII25.east) + (0,0.35)$)
							node [midway,above,yshift=7.5] {$r_{2}$}; 

 \end{tikzpicture}}
    \vspace*{-4mm}
    \caption{Coupled mass-spring-damper chain, consisting of two chains with $n_1$ and $n_2$ masses, respectively.}
    \label{fig:coupled_MSD-chain}
\end{figure}

\begin{minipage}[t]{0.475\textwidth}
\begin{figure}[H]
    \centering
%
%
\begin{tikzpicture}

\begin{axis}[%
width=2.25in,
height=1.5in,
at={(0.65in,0.4in)},
scale only axis,
xmode=log,
xmin=2000000,
xmax=100000000,
xminorticks=true,
xlabel style={font=\color{white!15!black}},
xlabel={total number of operations},
ymode=log,
ymin=1e-06,
ymax=0.1,
yminorticks=true,
ylabel style={font=\color{white!15!black}},
ylabel={global error at $t_{\mathrm{end}} = 2$},
axis background/.style={fill=white},
legend style={legend cell align=left, align=left, draw=white!15!black}
]
\addplot [color=blue, mark=triangle, mark options={solid, blue}]
  table[row sep=crcr]{%
2581133.33333333	0.000537617824238258\\
5149325.33333333	0.000134373232465437\\
10285709.3333333	3.35913561788482e-05\\
20558477.3333333	8.39771703253692e-06\\
41104013.3333333	2.09942160358329e-06\\
};
\addlegendentry{Strang splitting}

\addplot [color=red, mark=o, mark options={solid, red}]
  table[row sep=crcr]{%
5161683.33333333	0.00212805618433164\\
10281683.3333333	0.00053151015544051\\
20521683.3333333	0.000132846015933926\\
41001683.3333333	3.32095333563289e-05\\
81961683.3333333	8.30226017128028e-06\\
};
\addlegendentry{implicit midpoint rule}

\addplot [color=black, dashed, forget plot]
  table[row sep=crcr]{%
10240000	0.001192092895508\\
20480000	2.980232238769530e-04\\
};

\node[] at (axis cs: 2e+7,1.1e-3) {$p=2$};

\end{axis}

\begin{axis}[%
width=3in,
height=2in,
at={(0in,0in)},
scale only axis,
xmin=0,
xmax=1,
ymin=0,
ymax=1,
axis line style={draw=none},
ticks=none,
axis x line*=bottom,
axis y line*=left,
legend style={legend cell align=left, align=left, draw=white!15!black}
]
\end{axis}
\end{tikzpicture}%
    \caption{Coupled mass-spring-damper chain. Global error at $t_{\mathrm{end}} = 2$ vs.\ the total number of operations for the implicit midpoint rule ($\circ$) and the Strang splitting ($\triangle$).}
    \label{fig:MSD_chain_singlerate}
\end{figure}
\end{minipage} 
\hfill 
\begin{minipage}[t]{0.475\textwidth}
\begin{figure}[H]
    \centering
%
%
\begin{tikzpicture}

\begin{axis}[%
width=2.25in,
height=1.5in,
at={(0.65in,0.4in)},
scale only axis,
xmode=log,
xmin=4500000,
xmax=190000000,
xminorticks=true,
xlabel style={font=\color{white!15!black}},
xlabel={total number of operations},
ymode=log,
ymin=1e-07,
ymax=0.2,
yminorticks=true,
ylabel style={font=\color{white!15!black}},
ylabel={global error at $t_{\mathrm{end}} = 2$},
axis background/.style={fill=white},
legend style={legend cell align=left, align=left, draw=white!15!black}
]
\addplot [color=blue, mark=triangle, mark options={solid, blue}]
  table[row sep=crcr]{%
10793726.6666667	9.80234845749832e-05\\
21545726.6666667	2.45257815076803e-05\\
43049726.6666667	6.13268983702911e-06\\
86057726.6666667	1.53325052659882e-06\\
172073726.666667	3.83314825551436e-07\\
};
\addlegendentry{impulse method, m=10}

\addplot [color=red, mark=o, mark options={solid, red}]
  table[row sep=crcr]{%
5161683.33333333	0.00873248585642608\\
10281683.3333333	0.00238426570135564\\
20521683.3333333	0.000608553864477967\\
41001683.3333333	0.000152916929958133\\
81961683.3333333	3.82778526477199e-05\\
};
\addlegendentry{implicit midpoint rule}

\addplot [color=black, dashed, forget plot]
  table[row sep=crcr]{%
20480000	1.192092895507810e-05\\
40960000	2.980232238769530e-06\\
};
\node[] at (axis cs: 2.1e+7,2e-6) {$p=2$};

\end{axis}

\begin{axis}[%
width=3in,
height=2in,
at={(0in,0in)},
scale only axis,
xmin=0,
xmax=1,
ymin=0,
ymax=1,
axis line style={draw=none},
ticks=none,
axis x line*=bottom,
axis y line*=left,
legend style={legend cell align=left, align=left, draw=white!15!black}
]
\end{axis}
\end{tikzpicture}%
    \caption{Coupled mass-spring-damper chain. Global error at $t_{\mathrm{end}} = 2$ vs.\ the total number of operations for the implicit midpoint rule ($\circ$) and the impulse method with multirate factor $m=10$ ($\triangle$).}
    \label{fig:MSD_chain_multirate}
\end{figure}
\end{minipage}

\subsection{The singlerate case}
In a first simulation, we choose $n_1 = n_2 = 25$, $K_1 = K_2 = 50$, $K_{co} = 50$, $m_1 = m_2 = 0.3$, $r_1 = r_2 = 0.1$ and as initial value at $t=0$ we choose zero, apart from $x_{6}(0) = q_{1,3}(0) = 0.1$.
We perform numerical simulations on the time interval $t \in [0,2]$ using the implicit midpoint rule without any splitting approach \eqref{eq:implicit-midpoint}, as well as the Strang splitting \eqref{eq:Strang-splitting} with splitting \eqref{eq:diag_splitting} and approximation of the flows via the implicit midpoint rule. 
For both approaches, we investigate the global error at $t_{\mathrm{end}} = 2$ vs.\ the total number of operations  for different step sizes $h = 2^{-k},\ k=7,\ldots,11$. 
The results are depicted in Fig.\ \ref{fig:MSD_chain_singlerate}, emphasizing the efficiency of the splitting method. 

\subsection{The multirate case}
In a second example, we choose $n_1 = 5$, $n_2 = 45$, $K_{1} = 100$, $K_2 = 10$, $K_{co} = 10$, $m_1 = 0.1$, $m_2 = 0.4$, and $r_1 = r_2 = 0.1$. 
Consequently, the masses of the first chain are characterized by faster dynamics than the masses of the second chain, introducing multirate potential that motivates the use of multiple time stepping techniques like the impulse method \eqref{eq:impulse_method}.
Again, we choose zero as the initial value, apart from $x_6(0) = q_{1,3}(0) = 0.1$ and we perform numerical simulations on the time interval $t \in [0,2]$ using the implicit midpoint rule \eqref{eq:implicit-midpoint}, as well as the impulse method \eqref{eq:impulse_method} with splitting \eqref{eq:multirate-splitting} and approximation of the flows via the implicit midpoint rule. 
Analogously to the previous example, we compare the efficiency of the splitting approach compared to a straight-forward application of the implicit midpoint rule. 
As multirate factor, we choose $m=10$ and perform numerical simulations using the implicit midpoint rule \eqref{eq:implicit-midpoint} and the impulse method \eqref{eq:impulse_method} with (macro) step sizes $h^k$, $k=9,\ldots,13$. 
The results in Fig.\ \ref{fig:MSD_chain_multirate} highlight that the multiple time stepping approach results in a more efficient computational process.

\section{Conclusion and Outlook}\label{sec:conclusion_and_outlook}
\noindent In this work, we developed splitting methods tailored to coupled linear port-Hamiltonian systems.
The resulting algorithms are dissipativity-preserving and convergent of second order. 
The investigations also include multiple time stepping techniques, enabling the use of different step sizes for the subsystems.
Numerical results emphasize the efficiency of the proposed splitting methods compared to a straight-forward application of the implicit midpoint rule.

Future research will deal with the derivation of more advanced decomposition algorithms for linear PHS, particularly by using the (constant) Hessian $\nabla^2 H(x)$ of the Hamiltonian \cite{mönch2024commutatorbased}, enabling the derivation of fourth-order decomposition algorithms without using negative weights \cite{omelyan2003symplectic,schäfers2024hessianfree} and thus satisfying the energy balance at a discrete level.
We will moreover investigate the proposed splitting methods within the generalized additive Runge--Kutta (GARK) framework \cite{günther2023symplectic,schafers2023symplectic} that allows to couple certain steps which may result, for example, in improved stability properties.

\section*{Acknowledgements}
\noindent This work is supported by the German Research Foundation (DFG) research unit FOR5269 "Future methods for studying confined gluons in QCD".

\bibliographystyle{elsarticle-num} 
\bibliography{bibliography}
\end{document}